\documentclass[sn-mathphys-num] {sn-jnl} 

\usepackage{graphicx}%
\usepackage{multirow}%
\usepackage{amsmath,amssymb,amsfonts}%
\usepackage{amsthm}%
\usepackage{mathrsfs}%
\usepackage[title]{appendix}%
\usepackage{xcolor}%
\usepackage{textcomp}%
\usepackage{manyfoot}%
\usepackage{booktabs}%
\usepackage{algorithm}%
\usepackage{algorithmicx}%
\usepackage{algpseudocode}%
\usepackage{listings}%

\begin{document}

\title[Article Title] {
	When to use simulated annealing for solving CVRP? A case study of fuel deliveries in Poland
	}

\author*[1]
	{ \fnm {Vitalii} \sur {Naumov} }
	\email {vitalii.naumov@pk.edu.pl}
	\affil*[1] {
		\orgdiv {Transportation Systems Department},
		\orgname {Cracow University of Technology},
		\orgaddress {
			\street {Warszawska 24},
			\city {Krakow},
			\postcode {31155},
			\country {Poland}
		}
	}

\abstract{
	The paper addresses Capacitated Vehicle Routing Problem (CVRP) in the context of fuel delivery to gas stations. The CVRP aims to minimize total travel distance for a fleet with limited capacity. Fuel delivery, however, introduces unique complexities within the CVRP framework. We propose a novel approach that integrates the Simulated Annealing (SA) algorithm with a customized CVRP model specifically designed for gas station networks. This model incorporates real-world constraints like vehicle capacity, fuel demands at each station, and road network distances. The paper outlines the design of SA-based CVRP model for fuel delivery. We detail the objective function (minimizing distance) and the SA's exploration mechanism for generating candidate solutions. To assess its effectiveness, the proposed approach undergoes computational tests in Poland's gas station network serviced by the Samat transportation company. We compare the performance of our SA-based CVRP model with the conventional Mixed Integer Programming model for CVRP powered by Gurobi. The results aim to demonstrate the efficacy of the proposed SA-based heuristic in finding efficient routes for fuel deliveries.
	}


\maketitle

\section{Introduction}\label{sec1}

Despite the increasing popularity of electric and hybrid vehicles, traditional petrol-based fuels continue to dominate the road transport sector. This dominance is expected to persist for the foreseeable future, primarily due to factors such as the widespread availability of petrol and diesel infrastructure, the relatively higher cost of electric vehicles, and the limited range of electric vehicles compared to traditional fuel-powered vehicles. As a result, the delivery of petrol and diesel to gas stations remains a crucial link in the supply chain for road transport. Gas stations serve as the primary point of contact for consumers to refuel their vehicles, ensuring the uninterrupted operation of the transportation network. The efficiency and reliability of these delivery operations are essential to maintaining the smooth functioning of the transportation system and supporting economic activity.

Recent research has focused on addressing the complexities and challenges associated with fuel distribution, such as the inherent volatility of the product, the need for specialized transportation infrastructure, and the stringent safety regulations. For instance, Xiong et al. \cite{bib31} introduced a two-echelon vehicle routing problem to optimize liquefied natural gas (LNG) delivery from production terminals to filling stations, considering the impact of boil-off losses. Bonino et al. \cite{bib32} proposed a mixed-integer linear programming (MILP) approach for production and distribution planning in large-scale industrial gas supply chains.
Other studies have explored the impact of uncertainties, such as weather-related disruptions and fluctuating demand, on fuel distribution. Edwards et al. \cite{bib33} developed a stochastic scheduling model to optimize fuel delivery strategies during hurricanes, while Ghiami et al. \cite{bib41} investigated a deteriorating inventory routing problem for inland LNG distribution.

A pivotal challenge within the fuel delivery supply chain is the vehicle routing problem (VRP). This complex optimization problem involves determining the most efficient routes for a fleet of fuel tankers to deliver fuel to a network of gas stations. The goal is to minimize costs, such as fuel consumption, driver labor, and vehicle maintenance, while ensuring timely and reliable deliveries.

In the context of smaller-scale distribution, researchers have applied various heuristic and metaheuristic algorithms to solve VRP problems for liquefied petroleum gas (LPG) and LNG. Baihaqi and Fitria \cite{bib34} utilized the sweep and nearest neighbor algorithms to optimize the distribution of subsidized LPG, while Alkan et al. \cite{bib37} employed capacitated VRP with time windows (CVRPTW) and heuristic approaches to planning LPG distribution in the Aegean region (Turkey).
The integration of artificial intelligence (AI) and optimization techniques has also shown promise in improving LNG distribution efficiency. Assahla et al. \cite{bib38} explored the use of AI to optimize gas station replenishment, and Aqidawati et al. \cite{bib39} applied VRP to determine optimal routes for fuel distribution in Indonesia.

This paper investigates the capacitated VRP (CVRP) with a specific focus on the practical challenges of real-time fuel delivery dispatching. By developing a computationally efficient approach, this research aims to enhance the operational efficiency of fuel distribution networks. The proposed method can help mitigate the risk of fuel shortages at gas stations by enabling rapid response to changes in demand patterns. The computational experiments provided in this study intend to show the advantage of the proposed heuristic (in terms of computational time) over the traditional tools to solve CVRP for fuel deliveries considered by Kaleta et al. \cite{bib43}.

The remainder of this paper is organized as follows. Section \ref{sec2} reviews existing literature on CVRP, focusing on specific tools used to solve the optimization problem. Section \ref{sec3} formally defines the considered CVRP problem, incorporating specific constraints relevant to fuel delivery operations. In Section \ref{sec4}, the novel heuristic algorithm is proposed to efficiently solve the formulated CVRP problem. Section \ref{sec5} presents a comprehensive experimental evaluation of the proposed algorithm on a real-world gas station network. Finally, Section \ref{sec6} summarizes the key findings of this study and discusses potential avenues for future research.

\section {Literature review} \label{sec2}

CVRP is a classic NP-hard combinatorial optimization problem that has been extensively studied in the fields of operations research, logistics, and transportation. It involves determining optimal routes for a fleet of vehicles with limited capacity to serve a set of customers located at different locations. The primary objective is to minimize the total distance traveled while satisfying all customer demands, however, derivative objectives are also used oftentimes -- delivery costs, total delivery time, etc.

Over the years, researchers have proposed various approaches to tackle the CVRP, including exact methods, metaheuristic algorithms, and hybrid techniques. Each approach has its strengths and weaknesses, and the choice of method depends on the specific problem characteristics and computational resources.

While exact methods, such as branch-and-bound and cutting-plane algorithms, can guarantee optimal solutions, their computational complexity limits their applicability to small-scale instances. Practical use of the exact methods is usually not possible in the real-world conditions of fuel deliveries within the network of gas stations. Instead, metaheuristic algorithms have emerged as popular techniques for solving large-scale CVRP instances. These algorithms explore the solution space efficiently, often achieving near-optimal solutions within reasonable computational time.

Heuristics the most frequently used to solve CVRP include genetic algorithms (GA) \cite{bib2} \cite{bib3} \cite{bib5} \cite{bib8} \cite{bib10} \cite{bib14} \cite{bib18}, simulated annealing (SA) \cite{bib3} \cite{bib16} \cite{bib17} \cite{bib21} \cite{bib22} \cite{bib23} \cite{bib24} \cite{bib29}, and tabu search \cite{bib4} \cite{bib20} \cite{bib25} \cite{bib27} \cite{bib28} \cite{bib30}. However, more exotic approaches are presented in recent publications: an ant colony optimization (ACO) that mimics the behavior of ants foraging for food, using pheromone trails to guide the search \cite{bib6}, a particle swarm optimization (PSO) inspired by the social behavior of birds flocking, using particles to explore the possible solutions \cite{bib9}, and meerkat clan algorithm (MCA) inspired by the social behavior of meerkats, using clans to explore the solution space \cite{bib12}. The research by Asin-Acha et al. \cite{bib1} pioneers the integration of machine learning techniques to automate the selection and configuration of CVRP algorithms, leading to improvements in solution quality and computational efficiency.

The paper by Arshad et al. \cite{bib8} provides a comprehensive review of genetic algorithms for CVRP, highlighting their effectiveness in solving complex optimization problems.
Meniz and Tiryaki \cite{bib2} demonstrated the practical application of GA to solve CVRP in a real-world setting, optimizing the delivery network for a bread production facility. The authors found a reasonable solution for 215 shops serviced by 22 vehicles with the objective of minimization of distance.

Combining multiple metaheuristic algorithms can often yield improved performance of a heuristic method. A hybrid approach that integrates GA, SA, and local search techniques has been presented by Poonpanit et al. \cite{bib3} as a particularly effective tool to solve CVRP. Sehta and Thakar  \cite{bib7} address the specific challenge of eccentric CVRP instances, proposing a novel solution approach that combines convex hull and sweep algorithms with GAs. The paper by Kakkar et al.  \cite{bib5} proposes a two-phase approach using a modified genetic algorithm to efficiently solve CVRP instances: the clustering phase and the GA phase work synergistically to obtain high-quality solutions. Yu et al. \cite{bib29} introduced the Hybrid Vehicle Routing Problem (HVRP), considering the integration of plug-in hybrid electric vehicles into the CVRP framework: their study proposes a simulated annealing algorithm with a restart strategy to solve the HVRP, addressing the challenges of energy consumption and charging infrastructure. Ahmed et al. \cite{bib6} have proposed an enhanced ant colony system algorithm that incorporates $k$-nearest neighbor and local search techniques to overcome the limitations of traditional ACO algorithms. The paper by Sajid et al. \cite{bib9} presents an algebraic particle swarm optimization algorithm combined with SA to address the discrete nature of CVRP: the proposed algorithm has demonstrated superior performance compared to traditional methods.

The important direction in research related to the use of GAs for CVRP relies on possible improvements of the classical approach. The paper of Poonpanit et al. \cite{bib3} introduces a GA with local search techniques to enhance the performance of CVRP solvers: by addressing the issue of premature convergence, the proposed method may lead to improved quality of the obtained solutions. Arshad et al. \cite{bib8} also introduce novel crossover and mutation operators within GA tailored for CVRP.

\section{Problem formulation}\label{sec3}

The problem of delivering fuel to gas stations has many features and constraints that have a business nature. Many features may be considered by collecting the relevant and updated input data, such as demanded volumes of fuel and time moments the fuel needs to be delivered at stations. Other features may be evaluated for the existing distribution network, such as the matrices of the shortest distances and travel times, the service times, and the parameters of the available fleet (the number of vehicles and their capacities). The use of the problem data organized through the corresponding constraints reflecting business needs to obtain the best possible value of the objective function constitutes the contents of this section.

\subsection{Input data}\label{subsec31}

Input parameters are considered deterministic for the problem of routing fuel deliveries. Some of these parameters are deterministic by nature given that available infrastructure cannot be changed (the total number of stations, the distance matrix). However, other input parameters depend on the data provided by clients and may vary in practice (fuel amounts needed or delivery time) or are stochastic by nature (travel times, service times). Such non-deterministic characteristics are considered constant by using their expected (most probable) values to solve the optimization problem.

For the presented CVRP formulation, the following data is defined as an input to solve the problem for fuel deliveries to gas stations:

\begin{itemize}
\item $G$ is the number of gas stations to service
\item $H$ is the planning horizon in days
\item $N_k$ is the number of routes designated on $k$-th day, $k=1 \dots H$
\item $\mathbf{D}=\|d_{ij}\|$ is a matrix of the shortest distances between gas stations
\item $\mathbf{T}=\|t_{ij}\|$ is a matrix of travel times between gas stations
\item $\mathbf{S}=\|s_{ki}\|$ is a service time matrix: the sum of the forecasted refill time for the $k$-th day for the tanks at the $i$-th station
\item $\mathbf{U}=\|u_{kj}\|$ is the arrival times matrix: the latest arrival time of a vehicle at $j$-th station on $k$-th day
\item $\mathbf{Q}=\|q_{ki}\|$ is the forecasted demand for gas at $i$-th station on $k$-th day
\end{itemize}

Note that the number of routes planned for the given day is restricted by the available vehicles, whereas the number of days is defined by the planning horizon $H$. So, the total number $N$ of routes to be planned is defined as $N = \sum_{k=1}^{H}N_k$. Consequently, the number of rows in the matrices $\mathbf{S}$, $\mathbf{U}$, and $\mathbf{Q}$ are equal to the planning horizon $H$, and the number of columns is equal to the total number of serviced stations $G$. The matrices $\mathbf{D}$ and $\mathbf{T}$ are square, with the size equal to $G$.

\subsection{Objective function}\label{subsec32}

Traditionally, most formulations for CVRP use the total distance covered by vehicles as the objective function. That is justified in most cases, as business requirements usually include the minimization of the delivery costs or the minimization of the environmental impact which are derivative of the covered distance. The same reasoning is used in this study, and the objective function $z$ is presented as follows:

\begin{equation}
	z =  \sum\limits_{n=1}^{N} \sum\limits_{i=0}^{G}\sum\limits_{j=1}^{G+1} d_{i,j} \cdot x_{n,i,j} \rightarrow \min
	\label{eq1}
\end{equation}

where 
$x_{n,i,j}$ is a binary decision variable showing if a trip from $i$-th to $j$-th gas station is planned for $n$-th route:

\begin{align}
x_{n,i,j} &= 
	 \left\{ 
		\begin{aligned} 
			&1,  \text{if a trip  is assigned} \\
			&0, \text{otherwise}
		\end{aligned}
	\right.
\end{align}

As the decision variables are binary, the presented objective function is typical for Integer Programming (IP) problems. However, the need to consider the latest arrival times for each of the gas stations causes the use of auxiliary decision variables $y_{nj}$ -- times of arrival at $j$-th station for $n$-th route. As times of arrival are real numbers, the problem should be classified as a Mixed Integer Programming (MIP) problem, although decision variables $y_{nj}$ are not presented in the objective function.

\subsection{Constraints}\label{subsec33}

The model's constraints are multifaceted, encompassing both the traditional elements that define the optimal shape of delivery routes and the critical business imperatives that dictate the timing of deliveries to prevent fuel shortages at gas stations and the capacity limitations of the available fleet. The formulated constraints are exclusively based on the input data and defined decision variables.

\begin{enumerate}[1.]

\item The constraints that define a route configuration include the following conditions:

\begin{itemize}

\item Vehicles start the trip at the fuel hub:

\begin{equation}
	\bigwedge_{n=1,\dots,N}  \sum\limits_{j=1}^{G+1}x_{n,0,j} = 1
\end{equation}

To initiate the fuel delivery process, a vehicle's capacity must be filled with fuel at the designated fuel hub. Consequently, the starting point of each delivery route must be situated at this central hub.

\item Vehicles finish the route at the fuel hub:

\begin{equation}
	\bigwedge_{n=1,\dots,N}  \sum\limits_{i=0}^{G}x_{n,i,G+1} = 1
\end{equation}

A vehicle needs to return to the hub to proceed with deliveries on the next route. In practice, if a route is the last for the day, a dispatcher may manually assign the last route point that is different from a fuel hub. However, this condition should be used for the routing procedure to ensure minimum empty mileage at the obtained routes.

\item After arriving at the $i$-th station, a vehicle must depart from it:

\begin{equation}
	\bigwedge_{n=1,\dots,N} \bigwedge_{i=1,\dots,G} \sum\limits_{j=1}^{G}x_{n,j,i} - \sum\limits_{j=1}^{G}x_{n,i,j} = 0
\end{equation}

This constraint is necessary to guarantee that the routing procedure only outputs fully completed routes, avoiding any partial or incomplete solutions.

\item A vehicle cannot go next to the same station:

\begin{equation}
	\bigwedge_{n=1,\dots,N}  \bigwedge_{i=0,\dots,G} x_{n,i,i} = 0
\end{equation}

This constraint is essential to preclude the occurrence of duplicate elements within the derived solutions. As the diagonal elements ($d_{i,i} = 0$) of the shortest distance matrix may induce the optimization routine to select corresponding decision variables $x_{n,i,i}$, this constraint explicitly prohibits such selections.

\end{itemize}

\item Constraints defining arrival times:

\begin{itemize}
\item Restrictions on vehicle arrival times at gas stations:
\begin{equation}
	\bigwedge_{n=1,\dots,N_k} \bigwedge_{i=0,\dots,G} \bigwedge_{j=1,\dots,G+1} y_{n,i} - y_{n,j} + s_{k,i} + t_{i,j} - (1 - x_{n,i, j}) \cdot M \le 0
\end{equation}

The constraint imposes a time window restriction on route construction, ensuring that the arrival time at a subsequent station is contingent upon the service completion time at the preceding station and the inter-station travel time. The proposed formulation of this constraint uses a large positive number $M$ ("big $M$") to guarantee that for the $i$-th and $j$-th stations included in $n$-th route ($x_{n,i,j} = 1$) the above-described condition is satisfied, whereas for other stations, not included in the route ($x_{n,i,j} = 0$), the inequality defined by the constraint will always be true.

\item The route start time is set to zero:
\begin{equation}
	\bigwedge_{n=1,\dots,N} y_{n,0} = 0
\end{equation}

The initial departure time of the vehicle from the hub, excluding any tanking operations, is set to zero. Consequently, the latest arrival times for each serviced gas station must be dynamically recalculated based on this initial departure time (and defined as differences between the actual value of the needed arrival time and the planned start time of the route).

\item Upper bounds for times of arrival at gas stations:
\begin{equation}
	\bigwedge_{n=1,\dots,N_k}	\bigwedge_{j=1,\dots,G+1}   y_{n,j} \le u_{k,j}
\end{equation}

By imposing this constraint, the problem is transformed into a CVRP with time windows (CVRP-TW). The upper bounds of the time windows for each gas station are defined by their respective latest allowable arrival times, while the lower bounds are set to zero, allowing for flexibility in early arrivals.

\end{itemize}

\item Capacity constraints:

\begin{itemize}
\item The total demand at the stations visited (except the last one, where the vehicle tank is completely emptied) is less than the vehicle capacity:
\begin{equation}
	\bigwedge_{n=1,\dots,N_k}  \sum\limits_{i=1}^{G} \sum\limits_{j=1}^{G} x_{n,i,j} \cdot q_{k,i}  \le C
\end{equation}

\item The vehicle is emptied of fuel upon arrival at the destination point:
\begin{equation}
	\bigwedge_{n=1,\dots,N_k}  \sum\limits_{i=1}^{G} \sum\limits_{j=1}^{G+1} x_{n,i,j} \cdot q_{k,i}  \ge C
\end{equation}
where $C$ is a vehicle's capacity.

\end{itemize}

A homogeneous fleet of tank trailers, each with a standardized capacity of 39,000 liters, is considered in the presented model. This assumption reflects the prevalent vehicle configuration utilized in practical fuel delivery operations to gas stations in EU countries.

\end{enumerate}

\section{Proposed approach to solve the problem}\label{sec4}

A heuristic algorithm, consisting of two primary stages, is devised to solve the problem at hand. The initial stage focuses on generating a feasible initial solution, while the subsequent stage seeks to enhance the quality of the obtained solution through SA. A formal algorithmic description is provided in Algorithm \ref{algo1}.

The input data required for executing the procedure consists of the set of gas stations $\mathbf{\Omega}$ to be serviced, along with two key parameters: the number of iterations $R$ employed to generate the initial set of delivery routes, and the final temperature $T_{end}$ for the simulated annealing process, with an initial temperature of 1.

\begin{algorithm}
	\caption {Calculate delivery routes} \label{algo1}
	\begin{algorithmic}[1]
		\Function{routing}{$\mathbf{\Omega}$, $R$, $T_{end}$}
		\State $s_{MC} \gets \textbf{initial}(\mathbf{\Omega})$ \Comment{obtain an initial solution by running MC simulations}
		\State $z \gets \text{distance}(s_{MC})$
		\For{$i=2 \dots R$}
			\State $s_i \gets \textbf{initial}(\mathbf{\Omega})$
			\State $z_i \gets \text{distance}(s_i)$
			\If{$z_i < z$}
				\State $s_{MC} \gets s_i$
				\State $z \gets z_i$
			\EndIf		
		\EndFor
		\State $s_{SA} \gets \emptyset$ \Comment{improve the initial solution by running SA}
		\For{\textbf{each} $route$ \textbf{in} $s_{MC}$}
			\State $s_{SA} \overset{+}{\leftarrow} \textbf{SA}(route, 1, T_{end})$ \Comment{add the improved route to $s_{SA} $}
		\EndFor
		\State \textbf{return} $s_{SA}$
		\EndFunction
	\end{algorithmic}
\end{algorithm}

In the initial stage, a Monte Carlo (MC) simulation is executed to generate a diverse set of $R$ feasible delivery routes, each satisfying the problem constraints. The solution with the minimum total distance, denoted as $s_{MC}$, is identified and retained as the result of this stage. The proposed logic of the \texttt{initial} method for generating solutions is presented in Algorithm \ref{algo2}.

The second stage involves the application of the SA heuristic to each route within $s_{MC}$. The \texttt{SA} procedure, detailed in Algorithm \ref{algo3}, systematically explores the solution space by accepting both improving and non-improving moves, with the acceptance probability governed by a temperature parameter. The resulting set of improved routes $s_{SA}$ is returned as the final output of the algorithm.

The \texttt{initial} procedure from the set of stations to service $\mathbf{\Omega}$ randomly selects stations and adds them to routes. In each iteration of the inner loop, stations are randomly selected and added to a currently developed route if the vehicle capacity $C$ will not be surpassed and the delivery time will not be greater than the time $u_j$ required at the corresponding gas station. Once a station is assigned to a route, it is removed from the set $\mathbf{\Omega}$ to prevent duplicate assignments. The process of station selection and assignment continues until the total volume of fuel to be delivered on the route exceeds the vehicle's capacity.

The method \texttt{time($\omega$)} in the Algorithm \ref{algo2} denotes a function that returns the end time of servicing the station $\omega$ if it's included in the current route: it is calculated as the sum of travel times and service times for all previous stations visited at the route including the candidate station $\omega$.

\begin{algorithm}
	\caption {Obtain an initial solution} \label{algo2}
	\begin{algorithmic}[1]
	\Function{initial}{$\mathbf{\Omega}$}
	\State $routes \gets \emptyset$
	\While{$|\mathbf{\Omega}| > 0$} \Comment{generate new routes until the set of stations is empty}
		\State $route \gets \emptyset$
		\State $volume \gets 0$
		\While{$volume < C$}
			\State $j \gets \lfloor \text{rand}(0, |\Omega|) \rfloor $ \Comment{randomly select a station from $\Omega$}
			\If{$\text{time}(\omega_j) \le u_j \textbf{ and } volume + q_j \le C$}  
				\State $route \overset{+}{\leftarrow} \omega_j$ \Comment{add the station to the route}
				\State $\mathbf{\Omega} \overset{-}{\leftarrow} \omega_j$ \Comment{remove the station from the set of stations}
				\State $volume \gets volume + q_j$
			\EndIf
		\EndWhile
		\State $routes \overset{+}{\leftarrow} route$ \Comment{add the obtained route to the set of routes}
	\EndWhile
	\State \textbf{return} $routes$
	\EndFunction
	\end{algorithmic}
\end{algorithm}

The proposed \texttt{SA} procedure employs a classical approach to the simulated annealing heuristic. The input delivery route (given as an ordered set of gas stations) serves as the initial state, and its associated energy $E$ is calculated as the total distance of the route. The main loop of the algorithm iteratively generates new $candidate$ solutions until the temperature $T$ reaches a predefined end temperature $T_{end}$. In each iteration, a new solution is proposed by applying the \texttt{generate} procedure (see Algorithm \ref{algo4}) to the current solution. The new solution is accepted as the current solution if it results in a reduction of the total route distance ($E_c < E$). However, to avoid becoming trapped in local optima, the algorithm may also accept a solution that increases the total distance with a certain probability. This probability is determined by the Boltzmann distribution and is influenced by the temperature $T$. At the end of each iteration, the temperature is decreased according to the \texttt{decrease} method. Upon reaching the end temperature, the final route obtained through the SA process is returned as the output of the algorithm.

\begin{algorithm}
	\caption {SA procedure} \label{algo3}
	\begin{algorithmic}[1]
		\Function{SA}{$route, T_0, T_{end}$}
		\State $T \gets T_0$
		\State $E \gets \text{distance}(route)$
		\While{$ T > T_{end}$}
			\State $candidate \gets \textbf{generate}(route)$
			\State $E_c \gets \text{distance}(candidate)$
			\If{$E_c < E$}
				\State $route \gets cadidate$
				\State $E \gets E_c$
			\Else
				\If{$\exp(-\frac{E_c-E}{T}) \ge \text{rand}(0, 1)$}
					\State $route \gets cadidate$
					\State $E \gets E_c$
				\EndIf
			\EndIf
		\State $T \gets \text{decrease}(T)$
		\EndWhile
		\State \textbf{return} $route$
		\EndFunction
	\end{algorithmic}
\end{algorithm}

The \texttt{generate} method, a crucial part of the proposed SA procedure, operates by randomly selecting two distinct points along the delivery route provided as the argument. The selected points should satisfy a key condition: they cannot be neighbors within the route (this assumption requires that the route size should be more than 3 points: $|route|>3$). The order of stations located between these selected points is reversed then. This rearrangement results in a new route that services the exact same set of stations, albeit in a different sequence.

\begin{algorithm}
	\caption {Generate a new state (an alternative solution)} \label{algo4}
	\begin{algorithmic}[1]
	\Function{generate}{$route$}
		\State $n \gets |route|$
		\If{$n \le 3$}
			\State \textbf{return} $route$ \Comment{return the initial route if it's too small}
		\EndIf
		\State $i \gets \lfloor \text{rand}(0, n) \rfloor$
		\State $j \gets \lfloor \text{rand}(0, n) \rfloor$
		\While{$ | i - j | < 2$} \Comment{ensure the selected stations are not neighbors}
			\State $j \gets \lfloor \text{rand}(0, n) \rfloor$
		\EndWhile
		\If{$i > j$} \Comment{ensure that $i < j$}
			\State $\text{swap}(i, j)$
		\EndIf
		\State $rs \gets route[i:j]$ \Comment{select the route stations between $i$ and $j$}
		\State $\text{reverse}(rs)$ \Comment{reverse the order of the selected stations}
		\State \textbf{return} $route[:i] + rs + route[j:]$ \Comment{return the route with the reversed fragment}
	\EndFunction
	\end{algorithmic}
\end{algorithm}

The \texttt{rand($\alpha$, $\beta$)} method presented in the pseudocodes selects a value of the random variable uniformly distributed in the range between $\alpha$ (inclusively) and $\beta$ (exclusively). The \texttt{distance($\rho$)} method denotes a procedure that returns a total distance for the delivery routes $\rho$ provided as the argument -- the routes are defined as the ordered collections of gas stations to service.

\section{Results of experimental studies}\label{sec5}

The solution methodologies presented in this paper were implemented in Python for computational evaluation. The optimization model from Section \ref{sec3} was encoded using the \texttt{mip} package (\url{https://www.python-mip.com}). Alternatively, the heuristic approach described in Section \ref{sec4} was also implemented in Python. The complete source code for both implementations is publicly available on GitHub (\url{https://github.com/naumovvs/MC-SA-routing}).

To assess the performance of the proposed heuristic, a computational experiment was conducted using real-world data provided by Samat, a Polish fuel transportation company (\url{ https://samat.pl/en}). The data encompasses a large network of gas stations across Poland, where Samat delivers fuel from 16 terminals to over 500 stations servicing major brands like BP, Shell, and Circle K.

For this network, the shortest distances between all stations (matrix $\mathbf{D}$) and the corresponding travel times (matrix $\mathbf{T}$) were calculated using the HERE SDK tools (\url{https://www.here.com/platform/here-sdk}).

To comprehensively evaluate the proposed CVRP heuristic's computational efficiency, 300 data instances were generated based on Samat's historical data. Each instance incorporates information about the fuel demand (matrix $\mathbf{Q}$) at each gas station, the estimated service times (matrix $\mathbf{S}$), and the latest preferred delivery time (matrix $\mathbf{U}$) specified by the clients.

The numeric results shown in this Section refer to one of the tested instances (its parameters are presented in Table \ref{tab1}); however, the regularities discovered for this instance were also observed in each of the rest data instances.

\begin{table}[h]
	\caption{Instance characteristics} \label{tab1}
	\begin{tabular}{lll}
		\toprule
		Parameter & Unit & Value  \\
		\midrule
		The number of gas stations to service & stations & 65   \\
		Total amount of fuel to deliver & liters & 891,300   \\
		Mean amount of fuel to deliver per request & liters & 5369  \\
		Mean service time per station & minutes & 54  \\
		Planning horizon & days & 3  \\
		\botrule
	\end{tabular}
\end{table}

\begin{figure}[H]
	\centering
	\includegraphics[width=0.9\textwidth] {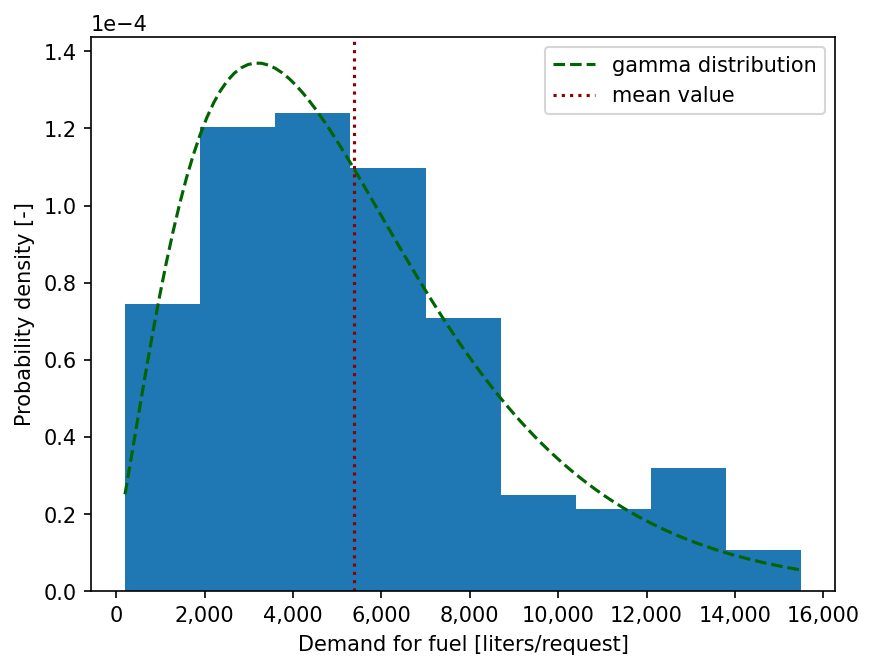}
	\caption {Distribution of the demand for fuel}
	\label {fig1}
\end{figure}

The distribution of fuel delivery requests per station is illustrated in Figure \ref{fig1}. The instance comprises 166 requests, ranging from 200 to 15,500 liters. As depicted in the histogram, the distribution of requested fuel amounts closely approximates a standard gamma distribution. This hypothesis is statistically supported by a chi-square Pearson's test with a significance level of 0.05.

Initially, the default CBC solver (COIN-OR Branch-and-Cut) within the \texttt{mip}-based code was employed to solve the optimization problem for the considered instance. However, this solver encountered difficulties and was unable to identify a feasible solution. To overcome this limitation, the Gurobi (\url{https://docs.gurobi.com}) solver was subsequently applied, successfully obtaining a feasible solution within a practical timeframe.

\begin{figure}[H]
	\centering
	\includegraphics[width=0.9\textwidth] {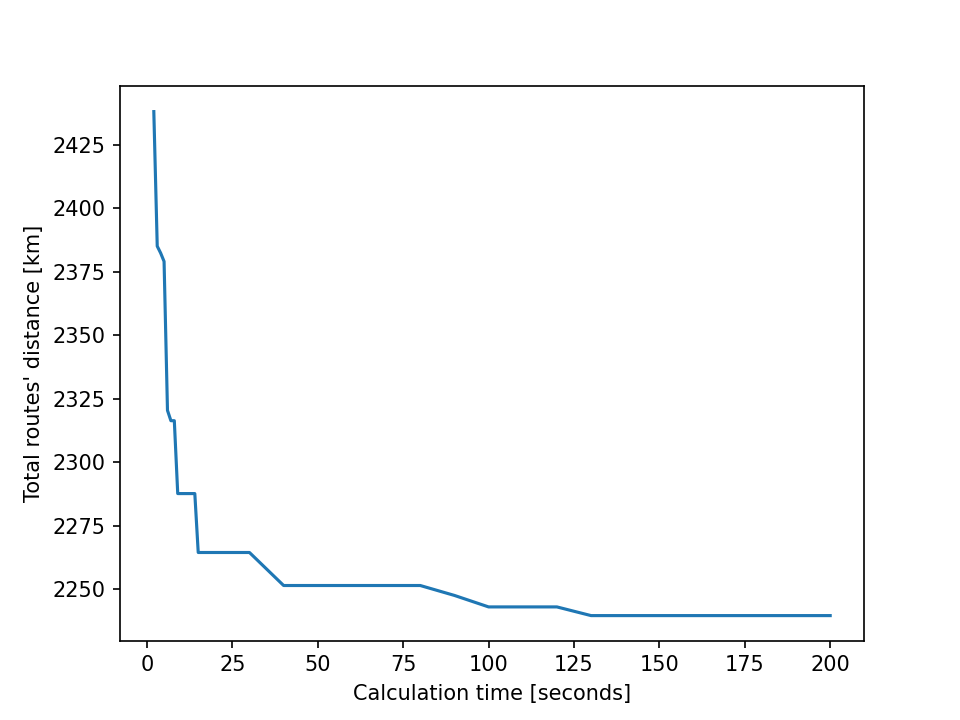}
	\caption {Optimization results using Gurobi}
	\label {fig2}
\end{figure}

Computational experiments were performed on a workstation equipped with an Intel Core i9 2.4 GHz processor and 16 GB of random-access memory. The optimization procedure, powered by Gurobi, was allowed to run for up to 200 seconds. Figure \ref{fig2} presents the resulting total route distances obtained within this time constraint.

As may be observed from the provided graphics, a substantial 7.12\% reduction in the total route distance, from an initial value of 2438 km to 2264 km, was achieved within the first 20 seconds of the Gurobi-powered optimization process. However, subsequent calculations yielded diminishing returns. Over the following 180 seconds, the total distance decreased by only an additional 1.10\% (from 2264 to 2239 km).

Considering the practical application of this decision-support tool in the context of transportation dispatching, where route planning is performed multiple times for varying demand levels, the ability to provide rapid solutions is crucial. The system must be capable of generating optimized routes within a few seconds, including the time required for data transmission and response reception. Given this practical constraint, the performance evaluation of the proposed SA procedure was focused on calculation times up to 20 seconds.

To assess the performance of the developed SA-based heuristic relative to the Gurobi-based optimization procedure, a comparative analysis was conducted for various numbers of Monte Carlo simulations. A larger number of simulations generally leads to a better initial solution but at the cost of increased computational time. As evident from Algorithm \ref{algo1}, the time complexity of the Monte Carlo phase exhibits a linear relationship with the number of trials, because of the same complexity of the loop iterations.

For the SA procedure, an initial temperature of 1 and a final temperature of 0.001 were employed. A cooling schedule based on inverse proportionality was adopted, with the temperature at the $i$-th iteration given by the rule $T_i = \frac{T_0}{i+1}$.

To ensure reliable performance evaluation of the developed SA-based heuristic, which incorporates stochastic elements from both the MC simulations and the subsequent SA phases, each experimental configuration was executed 100 times. The number of trials in the MC phase was varied as an input parameter to investigate its impact on solution quality and needed calculation time. For each experimental run, the computational time required for both the MC and SA phases was recorded. Additionally, the total routes' distance achieved at the conclusion of each phase was noted. To obtain statistically meaningful estimates, the average values of these metrics were calculated across the 100 repetitions for each configuration.

The selected results of the computational experiments for 3 levels of the number of trials in the MC phase ($R=1000$, $R=5000$, and $R=10,000$) are shown in Table \ref{tab2}.

Predictably, the mean duration of the SA stage doesn't change across different levels of $R$, as it depends on the number of serviced gas stations and the number of iterations in the SA main loop. Instead, the average duration of the MC stage can be approximated as about 1 second needed to generate 1000 initial solutions. As expected, the value of the objective function improves with the grown number of trials at the MC stage. It is noteworthy that the average total distance obtained after the SA stage with 1000 MC trials (2381 km) is nearly identical to the average objective function value achieved after the MC stage with 5000 trials (2386 km). However, the SA stage requires only 0.06 seconds of computation time, significantly less than the additional 4 seconds needed to generate 4000 more MC trials to achieve a comparable result.

\begin{table}[h]
	\caption{Results of the computational experiment} \label{tab2}
	\begin{tabular}{lllll}
		\toprule
		Parameter & Units & $R=1000$ & $R=5000$ & $R=10,000$ \\
		\midrule
		Mean total distance after the MC stage & km & 2446 & 2386 & 2365 \\
		Mean duration of the MC phase & seconds & 1.03 & 5.00 & 9.68  \\
		Mean total distance after the SA stage & km & 2381 & 2343 & 2326 \\
		Mean duration of the SA phase & seconds & 0.0581 & 0.0568 & 0.0577 \\
		\botrule
	\end{tabular}
\end{table}

The perspectives of practical use of the proposed heuristic should be estimated by comparing its computational efficiency with the results obtained by applying Gurobi as one of the world's leading optimization tools.

\begin{figure}[H]
	\includegraphics[width=0.9\textwidth] {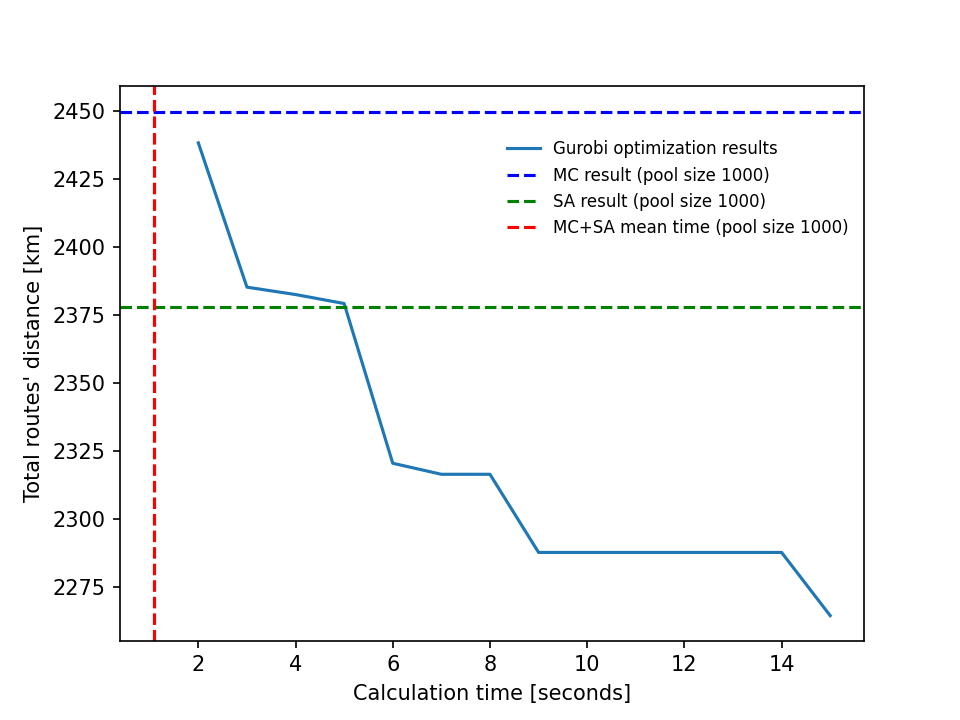}
	\caption {Results for 1000 try-outs in MC procedure}
	\label {fig3}
	\centering
	\includegraphics[width=0.9\textwidth] {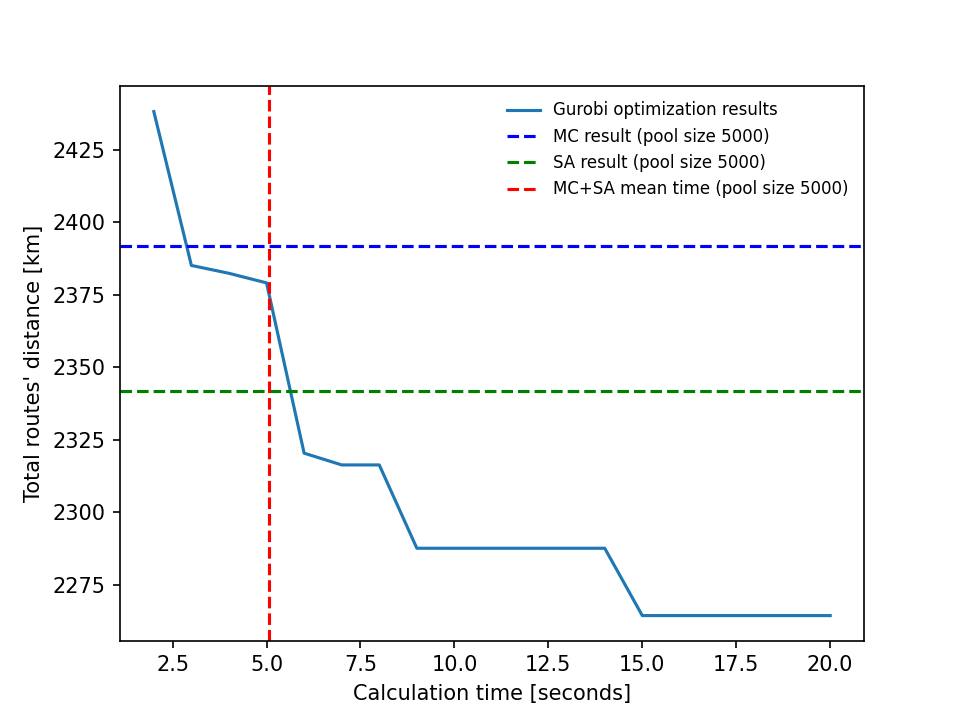}
	\caption {Results for 5000 try-outs in MC procedure}
	\label {fig4}
\end{figure}

The results for 1000 try-outs at the first stage presented in Table \ref{tab2} are also shown together with the Gurobi-based optimization results in Figure \ref{fig3}. The Gurobi solver initiated the optimization process and returned an initial feasible solution of 2438 km after a 2-second computation time. The proposed heuristic method, however, demonstrated superior performance by finding a solution of 2381 km in approximately 1 second. While the Gurobi-based optimization has the potential to converge to a better solution, it requires substantially longer computational time. As depicted in Figure \ref{fig3}, the Gurobi solver would necessitate 5 seconds of computation to attain the solution quality achieved by the heuristic method in a mere 1 second.

Figure \ref{fig4} illustrates the comparative performance of the heuristic and Gurobi-based approaches for 5000 Monte Carlo trials. The heuristic method, after its initial stage (approximately 5 seconds), provided a solution of 2386 km on average, which was slightly greater than the Gurobi-based solution of 2379 km. However, the subsequent application of the SA procedure within the heuristic framework significantly enhanced the solution quality to 2343 km. Although the Gurobi-based optimization eventually converged to a better solution of 2320 km after 6 seconds of computations, the developed heuristic method demonstrated a rapid initial improvement and a highly competitive final result. It is worth noting that running the developed routing method with about 5000 trials at the first stage represents the highest boundary of the heuristic effective use for the considered instance.

\begin{figure}[H]
	\includegraphics[width=0.9\textwidth] {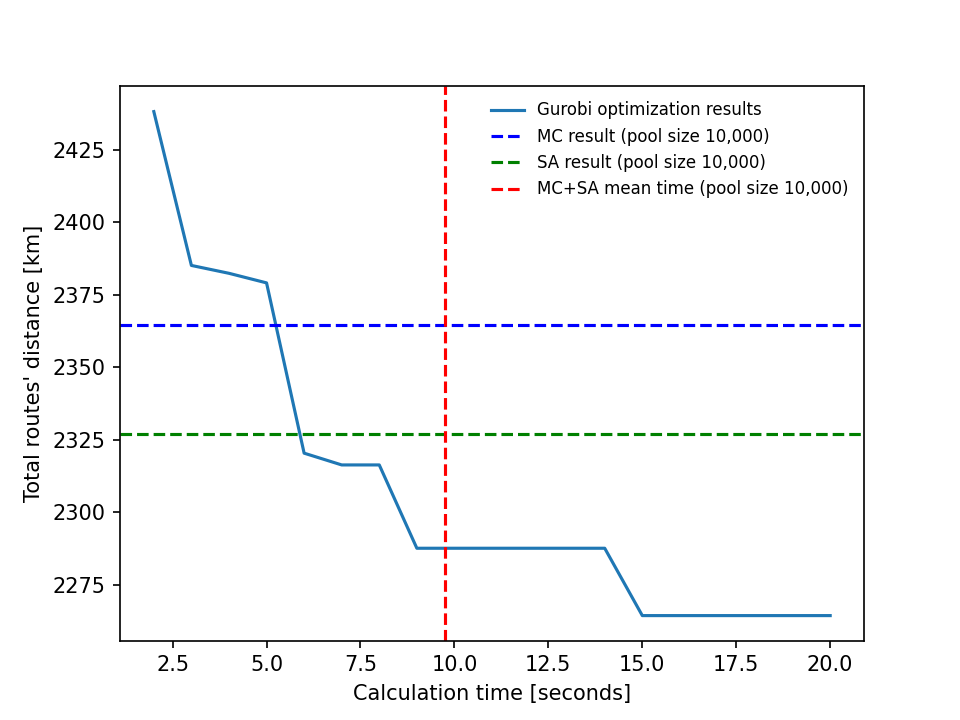}
	\caption {Results for 10,000 try-outs in MC procedure}
	\label {fig5}
\end{figure}

The performance of the developed heuristic evaluated through 10,000 Monte Carlo simulations and a comparison of the results obtained by the heuristic with those generated by a Gurobi-powered optimization model are presented in Figure \ref{fig5}. 

The figure illustrates that the proposed heuristic, on average, requires approximately 10 seconds to compute delivery routes with a total distance of 2326 km; while this is a respectable performance, the Gurobi-based optimization model outperforms the heuristic in terms of solution quality. After a similar computational time of 10 seconds, the Gurobi solver yields routes with a total distance of 2287 km, representing a 1.7\% improvement over the heuristic.

\section{Conclusion}\label{sec6}

This paper presents a novel approach to solving CVRP for fuel deliveries, by applying Monte Carlo simulations and the SA heuristic. The proposed method aims to determine a set of optimal vehicle routes that adhere to specific business constraints, such as vehicle capacity and delivery time bounds specified by clients.

In real-world logistics operations, dispatchers frequently need to re-plan routes in response to fluctuating demand levels. Therefore, a rapid and efficient routing solution is essential, so the developed SA-based heuristic is designed to generate high-quality routes within a matter of seconds.

Compared to traditional approaches that rely on powerful (and expensive for commercial use) optimization solvers like Gurobi, the proposed heuristic offers a compelling advantage: it can produce superior solutions in terms of total delivery distance, particularly within shorter computational timeframes. Computational experiments conducted on the Samat transportation company's data demonstrate that the SA-based heuristic consistently outperforms the Gurobi-powered model for up to 5 seconds of computation time. Beyond this threshold, the Gurobi solver may yield marginal improvements in total distance, suggesting that the heuristic provides a practical and effective solution.

The experimental results emphasize the practical value of the proposed heuristic. By striking a balance between solution quality and computational efficiency, this approach empowers dispatchers to make timely and informed routing decisions, ultimately leading to cost savings and operational optimization.

Based on the findings of this research, several promising directions for future investigation are identified:

\begin{enumerate}[1.]

\item
While the proposed heuristic has demonstrated effectiveness in the context of the Samat company's fuel delivery network, it is essential to assess its performance in other domains with potentially varying demand patterns and operational constraints. By conducting further experiments on diverse datasets, we can evaluate the heuristic's generalizability and identify potential limitations.

\item
To enhance the practical applicability of the proposed approach, it is worthwhile to explore the incorporation of additional business constraints. For instance, considering a heterogeneous vehicle fleet with varying capacities and fuel consumption rates can significantly impact route optimization. Prioritizing deliveries to specific gas stations can also be a crucial factor in optimizing operations. While the provided code repository includes an implementation for priority-based delivery, a comprehensive analysis of its impact on solution quality and computational efficiency is warranted.

\item
The SA stage plays a vital role in the proposed heuristic's performance. To further optimize the algorithm, a detailed investigation into the influence of different end temperatures and cooling schedules is necessary. By systematically varying these parameters, we can identify the optimal settings that lead to improved solution quality and computational efficiency. Additionally, exploring alternative neighborhood search strategies and acceptance criteria may yield further enhancements.

\end{enumerate}

\bmhead{Acknowledgements}
The study was co-funded by the National Center for Research and Development of Poland, under the R{\&}D project titled "Inventory routing optimization for petrol stations" (grant agreement no. POIR.01.01.01-00-0265/21).

\bibliography{sn-bibliography}

\end{document}